\documentclass[12pt]{amsart}
\usepackage{amscd,amssymb}
\usepackage[graph,frame,poly,arc]{xy}  
\usepackage[plainpages,backref,urlcolor=blue]{hyperref}

\theoremstyle{plain}
\newtheorem{thm}[subsection]{Theorem}

\newtheorem{cor}[subsection]{Corollary}

\theoremstyle{definition}
\newtheorem{rk}[subsection]{Remark}

\newtheorem{conj}[subsection]{Conjecture}

\numberwithin{equation}{section}
\newcommand{\V}{{\mathcal V}}
\newcommand{\f}{{\bf f}}
\newcommand{\g}{{\bf g}}
\newcommand{\dd}{{\bf d}}

\newcommand{\C}{\mathbb{C}}
\newcommand{\PP}{\mathbb{P}}

\begin{document}

	\title[Strong Lefschetz Property of complete intersections]{On Strong Lefschetz Property of 0-dimensional complete intersections and Veronese varieties}

        \author {Alexandru Dimca}
        \address{Universit\'e C\^ ote d'Azur, CNRS, LJAD, France and Simion Stoilow Institute of Mathematics,
        P.O. Box 1-764, RO-014700 Bucharest, Romania.}
       \email{dimca@unice.fr}

		\author{Giovanna Ilardi}
		\address{Dipartimento Matematica Ed Applicazioni “R. Caccioppoli” Universit`a Degli Studi	Di Napoli “Federico II” Via Cintia - Complesso Universitario Di Monte S. Angelo 80126 - Napoli - Italia}
		\email{giovanna.ilardi@unina.it}
		
		\author{Abbas Nasrollah Nejad}
		\address{Department of Mathematics, Institute for Advanced Studies in Basic Sciences (IASBS), Zanjan 45137-66731, Iran}
		\email{abbasnn@iasbs.ac.ir}
		


\begin{abstract} 
We show that the Strong Lefschetz Property in degree 1 for a homogeneous 0-dimensional complete intersection holds if the corresponding associated form, the Macaulay inverse systems, has a non-zero discriminant.
\end{abstract} 

\maketitle

\section{Statements of the main results} 

Let $S=\C[x_1, \ldots,x_n]$ be the polynomial ring with the usual grading and $n \geq 3$, $f_j \in S_{d_j}$ be a homogeneous polynomial of degree $d_j \geq 2$ for $j=1, \ldots,n$ such that the corresponding ideal
$$J(\f)=(f_1, \ldots ,f_n)$$
is a 0-dimensional complete intersection.
Consider the graded Artinian Gorenstein algebra $M(\f)=S/J(\f)$ with socle degree
\begin{equation}
\label{eq0} 
T=\sum_{j=1}^nd_j-n.
\end{equation}
In this paper we explore a possible new approach to a special case of the following well known conjecture, saying that the Artinian algebra $M(\f)$ has the Strong Lefschetz Property (SLP).

\begin{conj} \label{conj1}
For any 0-dimensional homogeneous complete intersection ideal $J(\f) \subset S$, any integer 
$k \in [0,T/2)$ and any generic linear form $\ell \in S_1$, the induced multiplication map
$$\ell^{T-2k}:M(\f)_k \to M(\f)_{T-k}$$
is an isomorphism.
\end{conj}

\begin{rk} \label{rk1}

(i) Let $f\in S_d$ be a reduced homogeneous polynomial of degree $d\geq 3$. By Euler formula the singular subscheme of the reduced hypersurface $X=V(f)\subseteq \PP^{n-1}$ is defined by the gradient ideal 
\[ J(f):=(f_1=\partial f/\partial x_1,\ldots, f_n=\partial f/\partial x_n).\]
The graded $\C$-algebra $M(f):=S/J(f)$ is called the Milnor algebra of $f$. If $f$ defines a smooth hypersurface in $\PP^{n-1}$, then $J(f)$ is generated by a regular sequence and hence $M(f)$ is an Artinian Gorenstein graded $\C$-algebra with socle degree $T=n(d-2)$. This setting is a special case of the case of 0-dimensional homogeneous complete intersections described above, and even in this case Conjecture \ref{conj1} is open in general.

\noindent (ii)
Conjecture \ref{conj1} (or a weaker form of it called the Weak Lefschetz Property) holds for the monomial ideals where
$$f_j=x_j^{d_j}  \text{ for all } j=1, \ldots,n,$$
 and for many other cases, including the case $n=2$ not considered here, see for instance \cite{BM,BFP,HMNW,I1,St,Wa3}. 
 
 \noindent (iii) Note that the case $k=0$ and any $n \geq 3$  holds for obvious reasons. Indeed, consider the degree $m$ Veronese embedding
\begin{equation}
\label{eq1} 
v_n^m:\PP(S_1) \to \PP(S_m), \  \   \ell \mapsto \ell^m.
\end{equation} 
Since the linear span of the Veronese variety 
$$\V^m_n=v_n^m(\PP(S_1))$$
is the whole projective space $\PP(S_m)$, this yields our claim if we take $m=T$.
However, even the cases $n=3$, or $k=1$, which  are among the simplest ones,  are open in general as far as we know. 
\end{rk}
In this note we address the case $k=1$ and give a sufficient condition on the ideal $J(\f)$ such that SLP holds in degree $k=1$ for the Artinian algebra $M(\f)$. Our setting uses  a new relation between 0-dimensional complete intersections, Veronese varieties and Macaulay inverse systems, which we recall briefly now, before stating our result.

For a 0-dimensional complete intersection $\f$ as above, one can associate a homogeneous polynomial $\mathrm{A}_\f$, called \textit{associate form} of $\f$, in the polynomial ring $R:=k[y_1,\ldots,y_n]$, where $y_i$ are dual variable to $x_i$~\cite{Alper-Isaev1,Alper-Isaev2,Eastwood-Isaev}. More precisely, denote by $\omega\colon \mathrm{Soc}(M(\f))\to  \C$ the linear isomorphism given by the condition $\omega(\overline{\mathrm{Jac}(\f)})=1$, where $\mathrm{Jac}(\f)$ is the determinant of the Jacobian matrix of the induced polynomial mapping $\f: \C^n \to \C^n$. The associate form $\mathrm{A}_\f$ is defined by the formula 
\begin{equation}
\label{eqAS1} 
\mathrm{A}_\f(y_1,\ldots,y_n)=\omega((y_1\overline{x_1}+\ldots+y_n\overline{x_n} )^T),
\end{equation}
where $\overline{x_i}\in M(\f)$ is the image of $x_i$. With the notation from Remark \ref{rk1} (i), it is clear that 
$$\mathrm{Jac}(\f)=\mathrm{Hess}(f),$$
is the determinant of the Hessian matrix of the polynomial $f$.

Given a homogeneous polynomial $F\in R$ of degree $T$, the \textit{apolar ideal} of $F$ is defined by 
\[ I_F=\mathrm{Ann}(F)=\{ g\in  S \ | \ g\circ F=0\}\subseteq R,\]
where the action of $S$ on $R$ is obtained by identifying $x_i$ with
$\partial /\partial y_i$.
Moreover, given a finitely generated graded Artinian Gorenstein $\C$-algebra $S/J$ with socle degree $T$, there is a homogeneous polynomial $F\in R_T$, unique up to scaling, such that $J=I_F$. Any such $F$
is called a \textit{Macaulay inverse system} for $S/J$. It is known that 
the associate form $A_\f(y_1,\ldots,y_n)\in R_T$ is a Macaulay inverse system for the Artinian algebra $M(\f)$, see \cite[Proposition 3.2]{Alper-Isaev1}. With the notation from Remark \ref{rk1} (i), we denote by
$$A_f(y_1,\ldots,y_n)\in R_T$$ the  Macaulay inverse system for the Milnor algebra $M(f)$.
With this notation our first main result is the following.

\begin{thm} \label{thmA3}
For any fixed multidegree $\dd=(d_1,\ldots,d_n)$ with
$$2 \leq d_1 \leq \ldots \leq d_n,$$
a generic
0-dimensional homogeneous complete intersection ideal $J(\f) \in S$ such that $\deg f_j=d_j$ for all $j=1,\ldots,n$ satisfies the following equivalent conditions.

\begin{enumerate}

\item The projective hypersurface defined by 
$$A_\f(y_1,\ldots,y_n)=0$$
is smooth.
\item
$\PP(J(\f)_{T-1}) \cap \V^{T-1}_n=\emptyset,$
where the intersection is taken in $\PP(S_{T-1}) $.

\end{enumerate}
\end{thm}
Note that the condition $(1)$ can be rephrased by saying that the discriminant $\Delta( A_\f(y_1,\ldots,y_n))$ of the polynomial $A_\f(y_1,\ldots,y_n)$ is non-zero, hence a numerical explicit condition.
With the notation from Remark \ref{rk1} (i), we have the following direct consequence of Theorem \ref{thmA3}.

\begin{cor} \label{corA3}
For any fixed degree $d \geq 3$,
a generic polynomial in $S_d$ satisfies the following equivalent conditions.

\begin{enumerate}

\item The projective hypersurface defined by 
$$A_f(y_1,\ldots,y_n)=0$$
is smooth.
\item
$\PP(J(f)_{T-1}) \cap \V^{T-1}_n=\emptyset,$
where the intersection is taken in $\PP(S_{T-1}) $.

\end{enumerate}
\end{cor}

The equivalence of the claims $(1)$ and $(2)$ in Theorem \ref{thmA3} follows from  \cite[Lemma 4.4]{Alper-Isaev1},
where the statement is for the Milnor algebra $M(f)$ as in our Corollary \ref{corA3} above, but the proof works in the general case, that is for the Artinian algebra $M(\f)$.
Moreover, Corollary  \ref{corA3} is proved in the same paper, see \cite[Proposition 4.3]{Alper-Isaev1}. To do this, the authors construct some explicit hypersurfaces with isolated singularities. By analogy, to prove Theorem \ref{thmA3}, we construct some 1-dimensional almost complete intersections, using some Bertini type theorems and generic choices of polynomials, see Theorem \ref{thmA1} below.

Our second main result, and motivation for this note, is the following, showing that for a specific  generic class of Artinian (resp. Milnor) algebras the SLP in degree 1 holds.
\begin{thm}
\label{thmSLP1}
The SLP holds in degree $k=1$ for any 0-dimensional homogeneous complete intersection $M(\f)$ for which the ideal $J(\f)$ satisfies the equivalent conditions (i) and (ii) in Theorem \ref{thmA3}.
In particular, the SLP holds in degree $k=1$ for any generic 0-dimensional homogeneous complete intersection  $M(\f)$.

\end{thm}
As a special case of this result, we get the following.

\begin{cor}
\label{corSLP1}
The SLP holds in degree $k=1$ for any Milnor algebra $M(f)$, when the polynomial $f$ of degree $d \geq 3$ is such the projective hypersurface $f=0$ is smooth and the Jacobian ideal $J(f)$  satisfies the equivalent conditions (i) and (ii) in Theorem \ref{thmA3}. In particular, the SLP holds in degree $k=1$ for any generic polynomial $f$.

\end{cor}

\section{Almost complete intersections and proofs of the main results}

Consider the ideal $K  \subset S$, generated by all the products
$x_ix_j$ for all $1 \leq i < j \leq n$. The following is a key technical result for our proofs.

\begin{thm} \label{thmA1}
For a given multidegree $\dd=(d_1,\ldots,d_n)$ such that
$$2 \leq d_1 \leq \ldots \leq d_n,$$
let $g_j $ be a generic polynomial in the vector space  $K_{d_j}$ of homogeneous polynomials in $K$ of degree $d_j$. 
Consider the ideal $J^j \subset S$ generated by the polynomials
$g_i$ for $1 \leq i \leq j$, where $1 \leq j \leq n$ and the corresponding
varieties $V^j=V(J^j) \subset \PP^{n-1}$. Then
$V^j$ is a smooth complete intersection of multidegree $(d_1,\ldots,d_j)$ and of dimension $n-1-j$ for any $j$ with  $1 \leq j\leq n-1$.
Moreover, $V^n=V(K)$ and hence consists of the $n$ points in $\PP^{n-1}$ given by the classes of the canonical basis
$e_1, \ldots, e_n$ of the vector space $\C^n$.
\end{thm}
\proof
We prove the first claim by induction on $j$. We start with the case $j=1$
and note that $V^1:g_1=0$ is a general member of the linear system defined by $K_{d_1}$. Using a Bertini type Theorem, see for instance \cite{GH}, p. 137, it follows that $V^1$ is smooth outside the base locus of $K_{d_1}$, which is precisely $V(K)$ since $K$ is generated in degree $2$. At a point in $V(K)$, say at the point $e_1=(1:0: \ldots :0)$ to fix the ideas, the linear part $\ell_1$ of $g_1$ localized at $e_1$ comes from the monomials $x_1^{d_1-1}x_2, \ldots, x_1^{d_1-1}x_n$. Since $g_1$ is generic, these monomials occur in $g_1$ with non-zero coefficients, and hence the corresponding linear part $\ell_1$ is non-zero. It follows that
$V^1$ is a smooth hypersurface of degree $d_1$. Assume that we have shown that $V^k$ is a smooth complete intersection for some $k<n-1$. To show that $V^{k+1}$ is a smooth complete intersection we proceed as follows. Let $p \in V^k$ be any point which is not in $V(K)$, and hence $g_{k+1}(p) \ne 0$, since $g_{k+1}$ is generic in $K_{d_{k+1}}$.
Since $V^k$ is a smooth complete intersection of dimension $>0$, it follows that $V^k$ is irreducible, and hence $V^{k+1}$ is a hypersurface in $V^k$, which is smooth at all the points outside $V(K)$ by the same argument as above. Now at a point in $V(K)$, let's say at the point $e_1$ to fix our ideas, the tangent space at $V^k$ is defined by
the vanishing of $k$ linear forms $\ell_1, \ldots, \ell_k$ in $x_2,\ldots,x_n$, which are obtained as we have explained in the case of $\ell_1$ above. The linear part $\ell_{k+1}$ of $g_{k+1}$ localized at $e_1$ comes from the monomials $x_1^{d_{k+1}-1}x_2, \ldots, x_1^{d_{k+1}-1}x_n$. Since $g_{k+1}$ is generic, these monomials occur in $g_{k+1}$ with generic coefficients, and hence the corresponding linear part $\ell_{k+1}$ is not in the linear span of the linear forms $\ell_1, \ldots, \ell_k$. This implies that $V^{k+1}$ is indeed a smooth complete intersection. This completes the proof of the first claim.

Now, to prove the second claim, note that $V^{n-1}$ is a 0-dimensional complete intersection, hence it is a finite set containing the set $V(K)$.
For any point $q \in V^{n-1} \setminus V(K)$, the condition $g_n(q)=0$ defines a hyperplane in the space $K_{d_n}$. Hence for a generic $g_n$ one has $g_n(q) \ne 0$ for any $q \in V^{n-1} \setminus V(K)$. This
completes the proof of Theorem \ref{thmA1}.

\endproof

It follows that the ideal $J({\bf g})=J^n \subset S$ generated by the polynomials
$g_1, \ldots, g_n$ is a 1-dimensional almost complete intersection, as  in \cite{DPop} and in \cite[Theorem 3.2]{Se}. Let $I({\bf g})$ be the saturation of the ideal $J({\bf g})$ with respect to the maximal ideal $(x_1, \ldots,x_n)$ and consider the graded $S$-module
$$N({\bf g})=\frac{I({\bf g})}{J({\bf g})}.$$
Then one has the following.

\begin{thm} \label{thmA2}
With the above notation,  the following holds.

\begin{enumerate}

\item $I(\g)=K$.

\item $N({\bf g})_1=N({\bf g})_{T-1}=0$.

\item $\dim (S/I({\bf g}))_{T-1}= \dim (S/J({\bf g}))_{T-1}=n$.

\end{enumerate}

\end{thm}
\proof
First we note that the ideal $I({\bf g})$ consists exactly of the polynomial in $S$ vanishing on the set $V(K)$, and hence $I({\bf g})=K$, see if necessary the description of the saturation of an ideal given in \cite[Section 2]{Bull13} using local vanishing conditions. In particular, $I({\bf g})_1=K_1=0$, which implies $N({\bf g})_1=0.$
To prove the rest of the second claim, we use the duality property of the graded $S$-module, see \cite[Theorem 3.2]{Se}, and get
$$\dim N({\bf g})_1=\dim N({\bf g})_{T-1},$$
and this completes the proof of the claim $(2)$. To prove the third claim, we recall that
$$\dim (S/I({\bf g}))_p = \deg I({\bf g})=n$$
for all $p \geq \sum_{j=1}^{n-1}d_j-(n-1)=T-(d_n-1)$, see \cite[Section 1, in particular the properties (P1)-(P4) and Lemma 1.1]{DPop}.
It follows that 
$$\dim (S/I({\bf g}))_{T-1}=n$$
and then the property $(1)$ implies that
$$\dim (S/J({\bf g}))_{T-1}=\dim (S/I({\bf g}))_{T-1}=n.$$
\endproof

\subsection{ Proof of Theorem \ref{thmA3}}

Note that the dimension one almost complete intersection ideal $J(\g)$
satisfies the condition 
$$\PP(J(\g)_{T-1}) \cap \V^{T-1}_n=\emptyset.$$
Indeed, Theorem \ref{thmA2} $(1)$ implies that
$$J(\g)_{T-1}=I(\g)_{T-1}=K_{T-1}.$$
Moreover, for any non-zero linear form $\ell \in S_1$, we clearly have
\begin{equation}
\label{C1} 
\ell^{T-1} \notin K_{T-1}.
\end{equation}
The condition  $\PP(J(\f)_{T-1}) \cap \V^{T-1}_n=\emptyset$ is an open condition, that is if it is satisfied for a 0-dimensional complete intersection $\f=(f_1, \ldots, f_n)$, then it is also satisfied for a nearby 0-dimensional complete intersection $\f '=(f_1', \ldots, f_n')$  with the same multidegrees  as $\f$. In particular, it follows that the claim in Theorem \ref{thmA3} holds for a given (multi)degree if it  holds for just one 0-dimensional complete intersection $\f=(f_1, \ldots, f_n)$ of that given multidegree.
Assume that this condition fails for all 0-dimensional complete intersections of multidegree $\dd$ and choose a sequence of such complete intersections $\f _m$ converging to the 1-dimensional complete intersection $\g$ constructed in Theorems \ref{thmA1} and \ref{thmA2}. Note that one has for the corresponding Artin Gorenstein algebras
$$\dim M(\f_m)_{T-1}=\dim M(\f_m)_1=n$$
for any $m$.  It follows that we get a sequence of linear subspaces $\PP(J(\f_m)_{T-1}) \subset \PP(S_{T-1})$ of codimension $n$. Passing to a subsequence if necessary, we may assume that this sequence converges to a codimendion $n$ subspace $E$ of $\PP(S_{T-1})$. Since
$\dim M(\g)_{T-1}=n$ as well by Theorem \ref{thmA2}, it follows that
$$E=J(\g)_{T-1}.$$
By our assumption, for any $m$ there is a point $p_m \in \PP(J(\f_m)_{T-1}) \cap \V^{T-1}_n$. By passing to a subsequence if necessary, we may assume that the sequence $p_m=(\ell_m')^{T-1}$ converges to a point $p=(\ell')^{T-1}$, which is both in $E=J(\g)_{T-1}$ and in $\V^{T-1}_n$.
This fact is in contradiction with \eqref{C1}, and in this way the proof of Theorem \ref{thmA3} is completed.

\endproof

\subsection{ Proof of Theorem \ref{thmSLP1}}
 Assume that SLP does not hold in degree $k=1$ for the Artinian algebra $M(\f)$. This means that for any generic linear form $\ell \in S_1$, there is another linear form
$\ell_1 \in S_1$ such that 
\begin{equation}
\label{eq2} 
\ell^{T-2}\ell_1 \in J(\f)_{T-1}.
\end{equation}
This is a closed condition, hence it would hold not only for generic $\ell$, but for all linear forms $\ell$. Recall that the projective tangent space to the Veronese variety $V^{T-1}_n$ at the point $p=\ell^{T-1}$ is given by
\begin{equation}
\label{eq3} 
T_pV^{T-1}_n= \{\ell^{T-2}\ell_1 \  | \ \ell_1 \in S_1\} \subset \PP(S_{T-1}),
\end{equation}
see for instance \cite[Section 1]{CGG}.
Consider now the linear projection of $\PP(S_{T-1})$ with center
$\PP(J(\f)_{T-1})$, which is a regular morphism
 \begin{equation}
\label{eq4} 
\pi:  \PP(S_{T-1}) \setminus \PP(J(\f)_{T-1}) \to \PP^{n-1}.
\end{equation}
Consider now the composition
$$\phi=\pi \circ v_n^{T-1}:\PP(S_1) \to \PP^{n-1}$$
which is again a regular morphism. The conditions \eqref{eq2} and \eqref{eq3} imply that the morphism $\phi$ is never a submersion, and hence Sard Theorem implies that $\phi$ is not a surjection.
This is a contradiction, because the regular morphism
$$\phi: \PP(S_1) =\PP^{n-1} \to \PP^{n-1}$$
must be given by some homogeneous polynomials
$$h_1, \ldots, h_n \in S_{T-1}$$
such that the corresponding ideal $I'=(h_1, \ldots, h_n )$ is a 0-dimensional complete intersection in $\C^n$ (otherwise $\phi$ is not defined everywhere on $\PP(S_1)$).
But such a 0-dimensional complete intersection $I'$ gives rise to surjective morphisms $\PP^{n-1} \to \PP^{n-1}$ as it is well known.

\newpage

\bigskip
\centerline  {\bf Conflict of Interests}

We declare that there is no conflict of interest regarding the publication of this paper.

\bigskip

\centerline {\bf Acknowledgments}

Alexandru Dimca was partially supported from the project ``Singularities and Applications'' - CF 132/31.07.2023 funded by the European Union - NextGenerationEU - through Romania's National Recovery and Resilience Plan.

Giovanna Ilardi was partially supported by GNSAGA-INDAM.

\end{document}